# Projective Chromatic Numbers

Adrian Rettich
Department of Mathematics
RPTU Kaiserslautern-Landau
Paul-Ehrlich-Straße 31
67663 Kaiserslautern, Germany
E-mail: rettich@rptu.de

Luke Serafin
Department of Mathematics
Cornell University
310 Malott Hall
14853 Ithaca, USA
E-mail: lss255@cornell.edu

January 2026

**Abstract**

We extend classical notions of definable colourability of graphs to the general projective setting and investigate whether known results, mainly about the $G_0$ dichotomy and the $2n+1$ conjecture, hold in the context of higher projective pointclasses. We establish that for $n \geq 2$, the presence of a $\mathbf{\Delta}^1_n$-definable well-order of the reals implies $\chi_{\mathbf{\Delta}^1_{n+1}}(G) = \chi(G)$ for all locally countable $\mathbf{\Delta}^1_n$-definable graphs $G$, and that the presence of a $\mathbf{\Delta}^1_2$-definable well-order of the reals implies $\chi_{\mathbf{\Delta}^1_2}(G) = \chi(G)$ for all locally countable Borel graphs $G$.

## 1 Introduction

Vertex colourings of graphs, that is, labellings where no two adjacent vertices can share a label, have been studied extensively both in the finite and the infinite case. Descriptive set theory has, ever since a seminal paper by Kechris, Solecki, and Todorčević [KST], combined this notion with its study of definable sets to spawn the new field of "descriptive combinatorics." In short, descriptive combinatorics looks for colourings where the labelling itself is not only a proper colouring, but satisfies additional assumptions, usually being a Borel-, Lebesgue-, or Baire-measurable map.

We investigate colourings definable by (boldface) projective pointclasses. For example, for a given graph one could look for the smallest cardinality of colours with which that graph can be coloured in a projectively defined way, or even at a specific rank of the projective hierarchy. This concept of chromatic number relative to an arbitrary level of the projective hierarchy has not to our knowledge been much investigated, and it significantly differs from much of classical descriptive combinatorics in that what is true can depend strongly on the model of set theory in which a statement is evaluated. For instance, the $\mathbf{\Delta}^1_2$ well-order of the reals which is available in $L$ can be used to construct a $\mathbf{\Delta}^1_2$ class $\aleph_0$-colouring of the graph $G_1$ on $2^\omega$ (which joins two vertices by an edge when they differ in finitely many coordinates), while in models of set theory where all sets of reals are Lebesgue-measurable (e.g. models of Martin's axiom for collections of dense sets of cardinality $\aleph_1$) the $\mathbf{\Delta}^1_2$ chromatic number of $G_1$ is uncountable because the Lebesgue-measurable chromatic number is.

In section 2 we introduce the classical notions before extending the definitions in a very natural way and collecting some initial observations about this expanded setting. Then in section 3 we investigate to what extent well-known results from descriptive combinatorics generalise to our case, including results on the so-called $2n+1$ conjecture and the $G_0$ dichotomy. Finally, in section 4 we present a first novel result about projective colourings, namely the fact that in the presence of a definable well-order of the reals at a certain level of the projective hierarchy some projective chromatic numbers of locally countable graphs collapse to the classical chromatic number.

We hope to establish with the present work a basis on which to investigate this aspect of descriptive combinatorics and present several open questions in section 5.

# 2 Preliminaries

**Definition 2.1.** We write $\mathbf{L}$ to denote the constructible universe and $\mathbf{V}$ to denote the von Neumann universe. If $A$ is a set, we write $\mathcal{P}(A)$ for the power set of $A$ and $|A|$ for the cardinality of $A$.

**Definition 2.2.** A topological space is called *Polish* if it is separable and completely metrisable.

The *Baire space*, written $\mathcal{N}$, is the Polish space $\omega^\omega$ equipped with the product topology, that is, the topology whose base consists of all sets of the form $\{ (w_i)_{i\in\omega} : \forall i < n : w_i = t_i \}$ for some finite sequence $(t_0, \dots, t_n)$ of length $n$ for some $n \in \omega$.

The Baire space is universal in the sense that it can be mapped continuously onto any Polish space.

**Definition 2.3.** Let $V$ be a topological space. A set $A \subseteq V$ is called *Borel* if it is an element of the $\sigma$-algebra generated by all open sets. We write $\mathcal{B}$ for this $\sigma$-algebra.

**Definition 2.4.** A subset of a Polish space is called *analytic* if it is the continuous image of a Polish space. We also call analytic sets $\mathbf{\Sigma}^1_1$ sets (note the boldface notation, which in this context is conceptually distinct from the lightface variant).

We now define the projective sets recursively.

**Definition 2.5.** A subset $A$ of a Polish space is $\mathbf{\Pi}^1_n$ if the complement of $A$ is $\mathbf{\Sigma}^1_n$.

A subset $A$ of a Polish space $V$ is $\mathbf{\Sigma}^1_{n+1}$ if there exists a set $A^\sharp \subseteq \mathcal{N} \times V$ such that $A^\sharp$ is $\mathbf{\Pi}^1_n$ and $A$ is the canonical projection of $A^\sharp$ onto $V$.

A subset $A$ of a Polish space is $\mathbf{\Delta}^1_n$ if it is both $\mathbf{\Sigma}^1_n$ and $\mathbf{\Pi}^1_n$.

For more on these classes, see for example [Kec]. The classes thus denoted (by boldface Greek letters) are known both as the (boldface) projective hierarchy and as the Lusin pointclasses (see [Mos]). They are not to be confused with the related *lightface* pointclasses, which can be thought of as an effective version of their boldface counterparts. An in-depth discussion is found in [Mos].

Note that due to the universality of the Baire space, it can be replaced by, for example, "there exists a Polish space $Y$ such that …" in the definition above.

**Definition 2.6.** A *graph* is a pair $G = (V, E)$, where $V$ is a set and $E \subseteq V \times V$ is a symmetric relation on $V$. The elements of $V$ are called the *vertices* of the graph, while the elements of $E$ are its *edges*. Two vertices $v, w$ are *adjacent* if and only if $(v, w) \in E$. The *open neighbourhood* $\mathrm{N}_G(v)$ of a vertex $v$ is the set of all vertices adjacent to $v$. The *closed neighbourhood* $\overline{\mathrm{N}}_G(v)$ of $v$ is its open neighbourhood plus $v$. In both cases, we drop the index $G$ if it is clear from context. For $W \subseteq V$ the *induced subgraph* on $W$ is $G[W] = (W, E \cap (W \times W))$. A *Polish graph* is a graph $(V, E)$ where $V$ is a Polish space. A *Borel graph* is a Polish graph $(V, E)$ where $E$ is a Borel subset of $V^2$.

Note the (non-standard) distinction between *Polish graphs* and *Borel graphs* which will allow us to state some results on graphs whose edge relation is not necessarily Borel.

**Definition 2.7.** Let $G = (V, E)$ be a graph, and let $A$ be a set of cardinality $\kappa$. A (proper) *$\kappa$-colouring* of $G$ is a function $V \to A$ such that no two adjacent vertices are assigned the same value (or "colour"). We call $G$ *$\kappa$-colourable* if such a colouring exists. The *chromatic number* of $G$, written $\chi(G)$, is the smallest cardinal $\kappa$ such that a $\kappa$-colouring of $G$ exists. Let now $c$ be a $\kappa$-colouring of a graph $(V, E)$. For $v \in V$, we call $\{ w \in V : c(w) = c(v) \}$ the *colour class* of $v$. Similarly, we call the set $\{ w \in V : c(w) = x \}$ the colour class of $x \in \kappa$. Given a Polish graph, a *Borel colouring* is a graph colouring where each colour class is a Borel set. The *Borel chromatic number* of a Polish graph $G$, written $\chi_{\mathcal{B}}(G)$, is the smallest cardinal $\kappa$ such that a Borel $\kappa$-colouring of $G$ exists. If a measure $\mu$ on the vertex set is given, then a *$\mu$-measurable colouring* is a graph colouring where each colour class is $\mu$-measurable, and the *$\mu$-measurable chromatic number* of a graph $G$, written $\chi_\mu(G)$, is the smallest cardinal $\kappa$ such that a $\mu$-measurable $\kappa$-colouring of $G$ exists. We write $\chi_\lambda(G)$ for the Lebesgue-measurable chromatic number of $G$.

Note that in their seminal paper on the topic, Kechris, Solecki and Todorčević [KST] define Borel colourings in a slightly different manner; in their original definition, a colouring $c : V \to X$ is Borel only if $X$ is a topological space as well (rather than an arbitrary set) and $c$ is a Borel map between these spaces. Our definition coincides with the more recent one used for example in [Ges] and has the advantage of distinguishing a richer universe of chromatic numbers in models where CH fails.[1] Geschke calls these *weak*

[1] Geschke's definition is strictly more discerning than the original definition, in which uncountable chromatic numbers always have value $\mathfrak{c}$.

*Borel chromatic numbers*, but we shall simply call them Borel chromatic numbers. Of course, a Borel colouring in the former sense is also a Borel colouring in the latter sense, but the reverse implication holds only in the countable setting.

While the classical chromatic number $\chi(G)$ is of course a lower bound for the Borel chromatic number, it is easy to see that the Borel chromatic number can be infinite even for acyclic graphs (see, for instance, the survey article [Ber]). Conley and Kechris [Con] showed that it is in fact possible to achieve an arbitrarily large gap between the two even on relatively simple graphs and even if both numbers are finite.

We extend these notions to the projective case in a natural way.

**Definition 2.8.** Let $\mathbf{P}$ be any Lusin pointclass. A $\mathbf{P}$ *class graph* is a graph $(V, E)$ where $V$ is a Polish space and $E \subseteq V \times V$ is in $\mathbf{P}$. Let $\kappa$ be a cardinal. A $\mathbf{P}$ *class $\kappa$-colouring* of a Polish graph $(V, E)$ is a $\kappa$-colouring where each colour class is in $\mathbf{P}$. Given a Polish graph $G$, its $\mathbf{P}$ *class chromatic number*, denoted $\chi_{\mathbf{P}}(G)$, is the smallest cardinal $\kappa$ such that a $\mathbf{P}$ class $\kappa$-colouring of $G$ exists. In all cases, we omit the word "class" when no confusion can arise, as in "$G$ is a $\mathbf{\Delta}^1_2$ graph".

If we want to emphasise that we are talking about $\chi(G)$ without any special property, we refer to it as the *classical chromatic number* of $G$.

Chromatic numbers (the classical as well as narrower notions) are well-defined for all graphs only if we assume the axiom of choice. For the remainder of this paper, we shall thus work in ZFC to minimise complications.

**Remark 2.9.** Because every Borel set is $\mathbf{\Delta}^1_2$ as well as Lebesgue measurable, we have for any Polish graph $G$ that $\chi_{\mathbf{\Delta}^1_2}(G) \leq \chi_{\mathcal{B}}(G)$ and $\chi_\lambda(G) \leq \chi_{\mathcal{B}}(G)$. However, it is consistent with ZFC that every $\mathbf{\Delta}^1_2$ set is Lebesgue measurable as well as that non-Lebesgue measurable $\mathbf{\Delta}^1_2$ sets exist [Mos].

We now show that in settings where the projective chromatic number is countable, we can restrict ourselves to talking about the pointclasses $\mathbf{\Delta}^1_n$.

**Lemma 2.10.** *Let $\mathbf{P} \in \{ \mathbf{\Sigma}, \mathbf{\Pi} \}$, and let $n \in \omega$, $n \geq 1$. Let $G$ be a Polish graph with $\chi_{\mathbf{P}^1_n}(G) \leq \aleph_0$. Then $\chi_{\mathbf{P}^1_n}(G) = \chi_{\mathbf{\Delta}^1_n}(G)$.*

*Proof.* Let $A$ be a colour class of a $\mathbf{P}^1_n$ class $\kappa$-colouring with $\kappa \leq \aleph_0$. Then $A$ is a $\mathbf{P}^1_n$ set. The complement of $A$ is the union of all other colour classes, of which there are only countably many and each of which is $\mathbf{P}^1_n$. Thus the complement of $A$ is $\mathbf{P}^1_n$ as well. □

**Lemma 2.11.** *Let $G$ be a Polish graph. Then $\chi_{\mathbf{\Sigma}^1_1}(G) = \chi_{\mathbf{\Pi}^1_1}(G) = \chi_{\mathcal{B}}(G)$.*

*Proof.* Let $G$ be any Polish graph. We show that $\chi_{\mathbf{\Sigma}^1_1}(G) = \chi_{\mathcal{B}}(G)$. The other claim is proven analogously. If $\chi_{\mathbf{\Sigma}^1_1}(G) \leq \aleph_0$, then this is lemma 2.10.

By Sierpiński's Theorem [Mos, theorem 2F.3], any $\mathbf{\Sigma}^1_2$ set of reals, and hence certainly any $\mathbf{\Sigma}^1_1$ set of reals, can be written as a union of at most $\aleph_1$ Borel sets. Consequently if $\chi_{\mathbf{\Sigma}^1_1}(G) > \aleph_0$, then

$$\chi_{\mathbf{\Sigma}^1_1}(G) \leq \chi_{\mathcal{B}}(G) \leq \aleph_1 \cdot \chi_{\mathbf{\Sigma}^1_1}(G) = \chi_{\mathbf{\Sigma}^1_1}(G). \qquad \square$$

It is worth noting that the above unlocks a specific technique for bounding Borel chromatic numbers: In order to bound the Borel chromatic number of a Polish graph from above, it suffices to find an analytic (or co-analytic) colouring, which may be easier than finding a Borel colouring.

**Lemma 2.12.** *Let $\mathbf{P} \in \{\mathbf{\Pi}^1_2, \mathbf{\Sigma}^1_2\}$, and let $G$ be a Polish graph with $\chi_{\mathbf{P}}(G) > \aleph_0$. Then $\chi_{\mathbf{P}}(G) = \chi_{\mathcal{B}}(G)$.*

*Proof.* The proof is the same as for lemma 2.11. $\square$

# 3 Generalisations of Classical Results

In this section, we consider some established facts about the Borel chromatic number and investigate whether or not they generalise to higher pointclasses.

## 3.1 Warm-Up Facts

The following fact is well-known for $\mathbf{P} = \mathbf{\Delta}^1_1$, but generalises trivially.

**Lemma 3.1.** *Let $\mathbf{P}$ be a Lusin pointclass and $G$ an acyclic $\mathbf{P}$ class graph with at least one edge. If there exists a $\mathbf{P}$-definable transversal of the connected components of $G$, then $\chi_{\mathbf{P}}(G) = 2$.*

We can also show that the bound known for degree-bounded Borel graphs (see [KST]) remains true higher up in the projective hierarchy.

**Theorem 3.2.** *Let $\mathbf{P}$ be a projective pointclass and $G$ a $\mathbf{P}$ class graph. If the degree of $G$ is bounded by $n \in \omega$, then $\chi_{\mathbf{P}}(G) \leq n + 1$.*

*Proof.* Our proof follows similar reasoning as [KST, chapter 4], though not all propositions shown there for Borel sets hold in our setting.

First, we notice that if $\mathbf{P}$ is a projective pointclass and $G = (V, E)$ is a locally countable $\mathbf{P}$ class graph with $\chi_{\mathbf{P}}(G) \leq \aleph_0$, then $G$ admits a maximal independent vertex set. Indeed: pick a minimal $\mathbf{P}$ class colouring of $G$ and call the colour classes $X_0, X_1, \ldots$. Inductively define $Y_0 = X_0$ and $Y_{n+1} = Y_n \cup (X_{n+1} \setminus \mathrm{N}(Y_n))$. Then clearly each $Y_n$ is independent, and thus so is $Y = \bigcup_n Y_n$, and for $v \in V \setminus Y$, there exists some $m > 0$ such that $v \in X_m$ (since $X_0 \subseteq Y$). Since $X_m \setminus \mathrm{N}(Y_{m-1}) \subseteq Y$, we must have $v \in \mathrm{N}(Y_{m-1}) \subseteq \mathrm{N}(Y)$, showing that $Y$ is maximal.

Let now $G = (V, E)$ be a $\mathbf{P}$ class graph with maximum degree $\Delta \in \omega$, and let $\{ U_n \}$ be a countable topological basis for $V$. Let $x \in V$. Since $\mathrm{N}(x)$ is finite, there is a least $c(x) \in \omega$ such that $x \in U_{c(x)}$ but $U_{c(x)} \cap \mathrm{N}(x) = \varnothing$. This is clearly a colouring, whence $\chi_{\mathbf{P}}(G) \leq \aleph_0$.

We can now use induction on $\Delta$. For $\Delta = 0$, there is nothing to do. Suppose thus that the claim has been shown for $\Delta - 1$. By what we have observed so far, we can find a maximal independent $\mathbf{P}$ set of vertices, say $Y$. Then by maximality of $Y$, the graph $G - Y$ has maximum degree $\Delta - 1$ and thus admits a $\mathbf{P}$ class $\Delta$-colouring. But since $Y$ is independent, we can complete this to a $(\Delta + 1)$-colouring of $G$. □

As in the general case, this bound is sharp, since any such graph can contain a complete graph on $n + 1$ vertices as a subgraph.

As in the Borel setting, graph homomorphisms can be a useful tool to bound chromatic numbers.

**Lemma 3.3.** *Let $G, G'$ be $\mathbf{P}$ class graphs and let $\varphi : G \to G'$ be a graph homomorphism which is also a $\mathbf{P}$ function. Then $\chi_{\mathbf{P}}(G) \leq \chi_{\mathbf{P}}(G')$.*

*Proof.* This is immediate from the definition of a graph homomorphism. □

## 3.2 Graphs Generated by Functions

**Definition 3.4.** Let $\mathbf{P}$ be a Lusin pointclass, let $\alpha \in \omega + 1$, let $V$ be a Polish space, and let $\{ f_i \}_{i \in \alpha}$ be a family of $\mathbf{P}$ functions from $V$ to $V$. We call the graph $(V, E)$ with

$$E := \{ (v, w) \in V \times V : v \neq w \wedge \exists i \in \alpha \, (f_i(v) = w \vee f_i(w) = v) \}$$

the graph *generated* by the family.

This definition coincides with established notions if $\mathbf{P}$ is taken to be $\boldsymbol{\Delta}^1_1$ (see for example [Pal]). The restriction to countable families of functions ensures that the resulting object is a $\mathbf{P}$ graph.

The following characterisation was shown very early in the study of Borel combinatorics.

**Theorem 3.5.** *A Polish graph is a locally countable Borel graph if and only if it is generated by a countable family of Borel functions* $\{ f_i \}$ *such that each* $f_i$ *is* $\leq \aleph_0$*-to-one, that is, for every* $v \in V$, $\left|f_i^{-1}(v)\right| \leq \aleph_0$.

*Proof.* See [KST, chapter 2]. □

Interestingly enough, this does not immediately generalise even to $\mathbf{\Pi}^1_2$ graphs, since the proof hinges on a certain uniformisation property[2].

At the time of writing, it is not known whether there exists a counterexample to theorem 3.5 in higher projective settings.

**Problem 3.6.** *Is there a model of set theory and a* $\mathbf{\Pi}^1_2$ *graph* $G$ *such that* $G$ *is locally countable, but not generated by any countable family of* $\mathbf{\Pi}^1_2$ *functions?*

The following proposition is included here as a first step towards finding such a counterexample.

**Lemma 3.7.** *Let* $E$ *be a relation generated by countably many* $\mathbf{P}$*-definable functions (in the sense of 3.4). Then* $E$ *is uniformised by a* $\mathbf{P}$ *set.*

*Proof.* Let $E$ be a relation generated by countably many $\mathbf{P}$-definable functions $f_0, f_1, \ldots$. Take

$$A_0 := \{ (x, f_0(x)) : f_0(x) \neq x \},$$

which is a $\mathbf{P}$ set. Define

$$A_{i+1} := \{ (x, f_{i+1}(x)) : f_{i+1}(x) \neq x \wedge \forall k \leq i (x = f_k(x)) \} \cup A_i,$$

another $\mathbf{P}$ set.

Then $\bigcup_{i\in\omega} A_i$ uniformises $A$ and is $\mathbf{P}$ as a countable union of $\mathbf{P}$ sets. □

Kanovei and Lyubetsky, in [KaV], showed that there exists a $\mathbf{\Pi}^1_2$ relation with countable sections which cannot be uniformised by any projectively definable set. This does not, however, resolve our question, since the relation constructed in [KaV] is not symmetric, and finding a symmetric relation with the same properties is an open problem.[3]

In fact, it is not even known whether the following is true.

[2]Namely, [KST] uses that any Borel relation with countable sections is generated by countably many Borel functions.

[3]After submission of this paper, Kanovei and Lyubetsky have actually constructed a symmetric relation with these properties, so we now know that the answer to problem 3.6 is “yes”.

**Problem 3.8.** *Let $G$ be a graph generated by finitely many $\mathbf{\Pi}^1_2$ functions. Is it then true that $\chi_{\mathbf{\Pi}^1_2}(G) \leq \aleph_0$?*

Indeed, even for graphs generated by a single function, the proof in the Borel case (see [KST]) hinges on the fact that given a Borel function on a Polish space, the topology can be refined into another Polish topology where the given function is continuous. No analogous result is likely to exist for higher projective pointclasses, although there is at present no counterexample to the above question.

In order to prove our main result, we use the following fact, which is slightly different from theorem 3.5.

**Lemma 3.9.** *Let $G = (V, E)$ be a locally countable Borel graph. Then there exists a countable family of Borel functions $\{ f_i : V \to V \}$ such that for every $(v, w) \in E$ there exists $i$ with $f_i(v) = w$.*

*Proof.* Let $G = (V, E)$ be a locally countable Borel graph. Then of course $E$ has countable sections, so by the Lusin-Novikov theorem there exist a countable family of Borel sets $\{ A_i \subseteq V \}$ and a countable family of Borel functions $\{ f_i : A_i \to V \}$ such that for every $(v, w) \in E$ there is $i$ with $v \in A_i$ and $f_i(v) = w$. The claim follows by noticing that each $f_i$ can be extended in a Borel way to have domain $V$ by having it be the identity outside of $A_i$ and that if $(v, w) \in E$ then also $(w, v) \in E$ because the edge relation is symmetric. $\square$

## 3.3 The $2n + 1$ Conjecture

The following is a well-known open problem in Borel combinatorics, called the $2n + 1$ conjecture.

**Conjecture 3.10.** *Let $G$ be a graph generated by $n \in \omega$ Borel functions. Then either $\chi_{\mathcal{B}}(G) = \aleph_0$, or $\chi_{\mathcal{B}}(G) \leq 2n + 1$.*

It is reasonable to conjecture the same for higher pointclasses, except we cannot necessarily expect the case of infinite chromatic number to be countable, as discussed in the previous section.

**Conjecture 3.11.** *Let $\mathbf{P}$ be a Lusin pointclass, and let $G$ be a graph generated by $n \in \omega$ functions, each of which is $\mathbf{P}$-definable. Then either $\chi_{\mathbf{P}}(G)$ is infinite, or $\chi_{\mathbf{P}}(G) \leq 2n + 1$.*

The Borel case has been shown for $n \in \{1, 2\}$, while the case $n = 3$ has been shown to have an upper bound of 8, one away from the conjectured bound of 7.

We show that the existing proofs can be lifted to the projective setting.

**Theorem 3.12.** *Let* $\mathbf{P}$ *be a Lusin pointclass, and let* $G$ *be a graph generated by a single* $\mathbf{P}$ *function. Then either* $\chi_{\mathbf{P}}(G) \in \{1, 2, 3\}$*, or* $\chi_{\mathbf{P}}(G)$ *is infinite.*

*Proof.* We can directly adapt the original proof from [KST]. We reproduce it here with slightly different wording and using arbitrary Lusin pointclasses rather than Borel sets.

Let $\mathbf{P}$ be a Lusin pointclass, let $f$ be a $\mathbf{P}$ function on a Polish space $V$, and let $G = (V, E)$ be generated by $f$.

For $n \in \omega$, $n \geq 2$, consider the graph $G_n$ on the space $n^\omega$ generated by the left-shift map $s_n : (x_0, x_1, \dots) \mapsto (x_1, x_2, \dots)$. For each $n$, this is a Borel function and thus also a $\mathbf{P}$ function.

We claim that for all $n \geq 2$, $\chi_{\mathcal{B}}(G_n) \leq 3$. For $n = 2$, we can partition $2^\omega$ into the Borel sets $\{ x \in 2^\omega : x \text{ starts with a finite odd number of zeroes} \}$, $\{ x \in 2^\omega : x \text{ starts with a finite odd number of ones} \}$, and the complement of the union of those two sets. Since each of those sets is independent in $G_2$, we have found a Borel 3-colouring.

We now proceed by induction. Suppose that we are given a Borel 3-colouring $c$ for $G_n$. Given $x \in (n+1)^\omega$, we explicitly assign to it a colour $c^+(x) \in \{0, 1, 2\}$. If $x$ contains infinitely many entries equal to $n$ and infinitely many entries not equal to $n$, we put

$$c^+(x) := \begin{cases} 0 & \text{if } x \text{ starts with a finite odd number of entries equal to } n \\ 1 & \text{if } x \text{ starts with a finite odd number of entries not equal to } n \\ 2 & \text{otherwise.} \end{cases}$$

The set of all such sequences clearly has no edges to the rest of the graph, and the three cases defined above are each independent.

Consider now sequences $x$ with only finitely many entries not equal to $n$. The set of all such sequences is a single connected component that is acyclic and countable, so it admits a Borel 2-colouring.

Finally, for every sequence $x$ with only finitely many entries equal to $n$ there exists $k_x \in \omega$ such that after the $k_x$th entry, $x$ contains no more entries equal to $n$. Set $c^+(x) := \left( c(s_{n+1}^{k_x}(x)) + k_x \right) \mod 3$. This is a Borel function, and two adjacent sequences are never assigned the same number by it: if $x$ has finitely many entries equal to $n$, then so do all of its neighbours.

If $x$ has finitely many but not zero entries equal to $n$, then $k_{s_{n+1}(x)} = k_x - 1$, but $s_{n+1}^{k_{s_{n+1}(x)}}(s_{n+1}(x)) = s_{n+1}^{k_x}(x)$, whence $x$ and $s_{n+1}x$ are not assigned the same colour. Finally, if $x$ contains no entries equal to $n$, then neither does $s_{n+1}(x)$, and by assumption, $c$ is a proper colouring of such vertices.

Thus, by induction, for every $n \in \omega$, we have $\chi_{\mathcal{B}}(G_n) \leq 3$.

Going back to our original graph $G$, by lemma 3.3 it suffices to show that either $\chi_{\mathbf{P}}(G)$ is infinite or there exists a $\mathbf{P}$-definable graph homomorphism from $G$ to some $G_n$.

If $\chi_{\mathbf{P}}(G)$ is finite, let $c$ be a $\mathbf{P}$ class $\chi_{\mathbf{P}}(G)$-colouring of $G$. Set $\varphi : G \to G_{\chi_{\mathbf{P}}(G)}, v \mapsto (c(v), c(f(v)), c(f^2(v)), c(f^3(v)), \dots)$. This is a $\mathbf{P}$ function, and if $v, w$ are adjacent in $G$, say $v = f(w)$, then $\varphi(v) = \varphi(f(w)) = s_{\chi_{\mathbf{P}}(G)}(\varphi(w))$. □

**Theorem 3.13.** *Let* $\mathbf{P}$ *be a Lusin pointclass, and let* $G$ *be a graph generated by two* $\mathbf{P}$ *functions. Then either* $\chi_{\mathbf{P}}(G) \leq 5$*, or* $\chi_{\mathbf{P}}(G)$ *is infinite.*

**Theorem 3.14.** *Let* $\mathbf{P}$ *be a Lusin pointclass, and let* $G$ *be a graph generated by three* $\mathbf{P}$ *functions. Then either* $\chi_{\mathbf{P}}(G) \leq 8$*, or* $\chi_{\mathbf{P}}(G)$ *is infinite.*

*Proof.* The proof for both of the preceding statements in the Borel setting can be found in [Pal]. Neither proof, while technically involved, uses the particular structure of the edge relation, and both proofs work just as well in our setting. □

## 3.4 The $G_0$ Dichotomy

In the context of Borel colourings, the graph $G_0$ plays a special role and actually fully characterises whether a locally countable graph has countable or uncountable Borel chromatic number. This can fail for graphs more complicated than Borel. In fact, [KST] shows that there is a model of ZFC where there exist two $\mathbf{\Pi}^1_1$ graphs $G$ and $H$ such that $\chi(G) > \aleph_0$ and $\chi(H) > \aleph_0$ but such that every graph that injects into both of them by graph homomorphisms must have countable $\mathbf{\Pi}^1_1$ chromatic number, whence of course it also has countable Borel chromatic number. Consequently, at least one of $G$ or $H$ is a graph with uncountable Borel chromatic number into which $G_0$ does not inject by *any* graph homomorphism, let alone a continuous one. We shall see that the $G_0$ dichotomy also strongly fails in the case of Borel graphs and $\mathbf{\Delta}^1_2$ colourings.

Recall the following definition.

**Definition 3.15.** We denote by $2^\omega$ the space of infinite binary sequences. Unless stated otherwise, we endow this space with the usual topology on sequence spaces (basic open sets are maximal collections of sequences which agree on an arbitrary but fixed finite prefix). We denote by $2^{<\omega}$ the set of finite binary sequences. Let $a \in 2^{<\omega}$ and let $b$ be any binary sequence. We say that $a$ is a *prefix* of $b$ if $|a| \leq |b|$ and $\forall i \leq |a| \, (a_i = b_i)$. A set $s \subseteq 2^{<\omega}$ is a *dense thread* if for every $k \in \omega$ it contains exactly one sequence of length $k$ and for every $a \in 2^{<\omega}$ there exists $a' \in s$ such that $a$ is a prefix of $a'$. Given a dense thread $s$ and $k \in \omega$, we denote the unique sequence in $s$ of length $k$ as $s_k$. Let now $s \subset 2^{<\omega}$ be a dense thread. We set

$$\begin{aligned} E := \{(a,b) \in 2^\omega \times 2^\omega : \exists k \in \omega \, (&s_k \text{ is a prefix of } a \text{ and } b, \\ &a_{k+1} \neq b_{k+1}, \\ &\forall i > k+1 \, (a_i = b_i))\} \end{aligned}$$

and define $G_s$ as the graph $(2^\omega, E)$.

We write $G_0$ when we care not from which specific dense thread the graph $G_s$ arises, even though the graphs generated by different threads are not usually isomorphic.

The $G_0$ *dichotomy* is then the following statement.

**Theorem 3.16.** *Let $G$ be a locally countable Borel graph and let $s \subset 2^{<\omega}$ be a dense thread. Then $\chi_{\mathcal{B}}(G)$ is uncountable if and only if $G_s$ injects into $G$ by a continuous graph homomorphism.*

*Proof.* This was shown in [KST], theorem 6.3, even for $\mathbf{\Sigma}^1_1$ graphs. □

In this sense, $G_0$ can be seen as the "smallest" Borel graph with uncountable Borel chromatic number.

One direction of course trivially holds for any pointclass above $\mathbf{\Delta}^1_1$.

**Lemma 3.17.** *Let $G$ be a locally countable Borel graph, let $s \subset 2^{<\omega}$ be a dense thread, and let $\mathbf{P}$ be a Lusin pointclass. Then if $\chi_{\mathbf{P}}(G) > \aleph_0$, $G_s$ injects into $G$ by a continuous graph homomorphism.*

*Proof.* This follows immediately from the fact that all projective chromatic numbers are bounded from above by the Borel chromatic number. □

We show, however, using a result from section 4, that the $\mathbf{\Delta}^1_2$ chromatic number of $G_0$ can vary wildly across different models of ZFC.

**Theorem 3.18.** *Let $s \subset 2^{<\omega}$ be a dense thread. Then the following are true:*

- *There is a model of ZFC in which $\chi_{\mathbf{\Delta}^1_2}(G_s) = 2$.*
- *There is a model of ZFC in which $\chi_{\mathbf{\Delta}^1_2}(G_s) > \aleph_0$.*

*Proof.* For the first statement, observe that $G_s$ is locally countable, so by theorem 4.4, $\mathbf{L}$ is a model where $\chi_{\mathbf{\Delta}^1_2}(G_s) = 2$.

For the second statement, choose any model where all $\mathbf{\Delta}^1_2$ sets of reals are Baire-measurable and the covering number of the meagre ideal is at least $\aleph_2$ (see, for example, [Bar, section 2.4]). The proof in [KST] that the Borel chromatic number of $G_0$ is uncountable shows that it is in fact at least as large as the covering number of the meagre ideal, so at least $\aleph_2$ in our chosen model.

Suppose now its $\mathbf{\Delta}^1_2$ chromatic number were countable. Then by Sierpiński's Theorem [Mos, theorem 2F.3], given a $\mathbf{\Delta}^1_2$ colouring of size $\kappa \leq \aleph_0$, we could split each colour class into at most $\aleph_1$ Borel sets, yielding a Borel colouring of size $\aleph_1 \cdot \kappa \leq \aleph_1$, a contradiction to $\chi_{\mathcal{B}}(G_0) \geq \aleph_2$. □

# 4 Models with Projective Well-Orders

Our main goal in this section is to show the following.

**Theorem 4.4.** *If there exists a $\mathbf{\Delta}^1_2$ well-order of $\mathcal{N}$, then for every locally countable Borel graph $G$, $\chi_{\mathbf{\Delta}^1_2}(G) = \chi(G)$.*

The key step of the proof is to give a "nice" $\mathbf{\Sigma}^1_2$ well-order of the vertices of $G$: on each connected component, we want our well-order to induce a well-order of order type at most $\omega$. Once such an enumeration is available, a $\mathbf{\Sigma}^1_2$ colouring of $G$ is obtained by sorting a set of candidate functions and defining the colouring on each connected component by induction (lemma 4.2).

Note that a $\mathbf{\Sigma}^1_2$-definable linear order is necessarily $\mathbf{\Delta}^1_2$-definable as well, though we shall not need this fact because lemma 2.10 applies to all cases discussed in this section.

A nice well-order is obtained by ordering elements first by the least parameters, according to the order of construction, used to produce them by application of Gödel functions, and then by a well-order in type $\omega$ of these Gödel functions (lemma 4.3).

We first show that if we have a sufficiently nice well-order, then the chromatic numbers from that part of the projective hierarchy onwards collapse. This works on any level above 1, not only for $\mathbf{\Delta}^1_2$.

**Definition 4.1.** Let $(V, E)$ be a graph. For $v \in V$, we write $𐤇_v$ for its connected component, that is, the set of all vertices reachable from $v$ by a finite path.[4]

**Lemma 4.2.** *Let $n \geq 2$, let $G = (V, E)$ be a locally countable $\mathbf{\Delta}^1_n$ graph, and let $\preceq$ be a $\mathbf{\Sigma}^1_n$-definable well-order of $V$. If, for every $v \in V$, the order type of the connected component of $v$ ordered by the restriction of $\preceq$ is at most $\omega$, then $\chi_{\mathbf{\Sigma}^1_n}(G) = \chi(G)$.*

*Proof.* Let $n \geq 2$, let $G = (V, E)$ be a locally countable $\mathbf{\Delta}^1_n$ class graph, and let $\preceq$ be a $\mathbf{\Sigma}^1_n$-definable well-order on $V$ such that for every $v \in V$, the order type of the connected component of $v$ ordered by the restriction of $\preceq$ is at most $\omega$.

Define, for $v \in V$, $𐤊(v) := \operatorname{ot}\{\, w \in 𐤇_v : w \prec v \,\}$.[5] *Nota bene*, since the order type of $(𐤇_v, \preceq)$ is at most $\omega$, each $𐤊(x)$ is a natural number, and no two vertices in the same connected component are assigned the same number. Note also that for the same reason, the function $𐤊 : V \to \omega$ can be defined inductively in a $\mathbf{\Sigma}^1_n$ way.

We first observe that since each connected component is countable, the function $𐤊$ itself is a proper $\aleph_0$-colouring. Thus if $\chi(G)$ is infinite, the claim holds. From now on we assume without loss of generality that $\chi(G) < \aleph_0$.

We shall construct a $\mathbf{\Sigma}^1_n$-definable $\chi(G)$-labelling of $V$ and then show that this is a proper colouring.

We write $C := \chi(G)^\omega$, a Borel subset of $\mathcal{N}$.

Let $v \in V$. We call a $a \in C$ a *$v$-palette* if it has the following property:

$$\forall x \in 𐤇_v \, \forall n \in \omega \, (a(𐤊(x)) = n \to \forall y \in 𐤇_x \, (y \in \mathrm{N}(x) \to a(𐤊(y)) \neq n)).$$

For every $v \in V$, at least one $v$-palette exists: we have assumed without loss of generality that $\chi(G)$ is finite. Let $c : V \to \chi(G)$ be any classical $\chi(G)$-colouring of $G$, and let $v \in V$. Then $c' : \omega \to \chi(G)$ with $c'(n) = m$ if and only if

$$\exists x \in 𐤇_v \, (𐤊(x) = n \wedge c(x) = m) \vee m = 0 \wedge \neg \exists x \in 𐤇_v \, (𐤊(x) = n)$$

[4] 𐤇 is the Phoenician letter *het*.
[5] 𐤊 is the Phoenician letter *kaf*.

is a $v$-palette. Furthermore, because the order type of each $(\text{目}_v, \preceq)$ is countable, the palette condition is $\mathbf{\Sigma}^1_n$-definable, whence the set of all $v$-palettes is a $\mathbf{\Sigma}^1_n$ subset of $C$. We shall denote this set by $C_v$.

Note that if $v, w$ are connected, then $C_v = C_w$.

Let $v \in V$. We define a function $c(v)$ inductively as follows. Suppose $c$ has been defined on at least the set $\{\, x \in \text{目}_v : x \prec v \,\}$. Then we set $c(v)$ to be the smallest natural number $n$ such that there exists $d \in C_v$ with

$$d(v) = n \wedge \forall x \in \text{目}_v \, (x \prec v \rightarrow d(x) = c(x)).$$

By construction, this is well-defined and a $\mathbf{\Sigma}^1_n$ function (the only quantifier over an uncountable set being the existential one over $C_v$).

It is also a proper colouring: let $v, w \in V$ be adjacent. Then there exists $d \in C_v = C_w$ such that $c(v) = d(v)$ and $c(w) = d(w)$. Since $d$ is a palette, we must have $d(v) \neq d(w)$.

Because no palette can use more than $\chi(G)$ colours, this proves the claim. □

We must now show how to apply this lemma to our case. Mansfield's theorem [Jec, Theorem 25.39] limits the complexity of the models we must consider. Here, $\mathbf{L}[x]$ is the constructible universe relative to $x$ (see, for example, [KaA, section 1.3]). Due to this fact, to prove theorem 4.4 it suffices to work in a model of the form $\mathbf{L}[x]$ for some $x \in \mathcal{N}$.

**Lemma 4.3.** *Let $x \in \mathcal{N}$ such that $\mathcal{N} \subseteq \boldsymbol{L}[x]$, and let $G$ be a Borel graph. In the model $\boldsymbol{L}[x]$ there is a $\mathbf{\Sigma}^1_2$-definable well-order $\preceq$ of $\mathcal{N}$ with the property that $\preceq$ restricted to any one connected component of $G$ has order type at most $\omega$.*

*Proof.* Let $x \in \mathcal{N}$ such that $\mathcal{N} \subseteq \mathbf{L}[x]$, and let $G$ be a Borel graph. It is well-known (see, for example, [Jec, theorem 25.26] and the discussion following it) that in the model $\mathbf{L}[x]$, the order of construction is a $\Delta^1_2(x)$ well-order. We use the fact that a subset of $\mathcal{N}$ is $\mathbf{\Sigma}^1_2$ if and only if it is $\Sigma^{\mathrm{HC}}_1$ [Jec, lemma 25.25].

Since $G$ is a Borel graph, by lemma 3.9 there exists a countable sequence $F$ of Borel functions $V \rightarrow V$ such that for every $(v, w) \in E$ we can find $f$ in the sequnce $F$ with $f(v) = w$.

We define an alternative well-order on $\mathcal{N}$ as follows. For $x \in \mathcal{N}$, write $\bar{x}$ for the least (lexicographically in the order of construction) set of parameters defining $x$ via application of Gödel functions such that $F \in \bar{x}$. In other words, we use the usual order of parameters, but additionally always allow the countable sequence $F$ in the parameter list. This in particular has the effect of making the use of elements of $F$ in definitions "free."

We claim that if $x, y \in \mathcal{N}$ are in the same connected component, then $\bar{x} = \bar{y}$. Indeed: if $x \in \mathcal{N}$ is definable from the parameter set $\bar{x}$ and $y \in \mathcal{N}$ is adjacent to $x$, then as noted before, we can find $f \in F$ (definable from $\bar{x}$ since $F \in \bar{x}$) such that $f(x) = y$, and thus $y$ is definable from $\bar{x}$ as well. Since any two vertices in the same connected component must be connected by a finite sequence of adjacent vertices, the claim follows.

We can now order $\mathcal{N}$ as follows: for $x, y \in \mathcal{N}$, if $\bar{x}$ is less than $\bar{y}$ (lexicographically in the order of construction), then we say $x \mathrel{\bar{\prec}} y$. This is in $\Sigma_1^{\mathrm{HC}}$ since the original order of construction is. If $\bar{x} = \bar{y}$, then we order by the least Gödel functions defining $x, y$ from $\bar{x}$.

It is easy to see that this is a well-order; it is $\Sigma_1^{\mathrm{HC}}$ since the original order of construction is, and by the result cited earlier, it is $\boldsymbol{\Sigma}_2^1$.

This well-order restricted to any one connected component is ordered solely by Gödel functions, which we take to be computably enumerated in order type $\omega$, and thus the order type of the connected component is no more than $\omega$. □

We thus arrive at the main result.

**Theorem 4.4.** *If there exists a $\boldsymbol{\Delta}_2^1$ well-order of $\mathcal{N}$, then for every locally countable Borel graph $G$, $\chi_{\boldsymbol{\Delta}_2^1}(G) = \chi(G)$.*

*Proof.* Let $G = (V, E)$ be a locally countable Borel graph in such a model. By Mansfield's theorem [Jec, Theorem 25.39], in this model we have $\mathcal{N} \subseteq \mathbf{L}[x]$ for some $x \in \mathcal{N}$. Thus by lemma 4.3 the model contains a $\boldsymbol{\Sigma}_2^1$ well-order of $V$ which, restricted to any one connected component, enumerates the vertices of that connected component. By lemma 4.2, we get $\chi_{\boldsymbol{\Sigma}_2^1}(G) = \chi(G)$. Since $G$ is locally countable, we have $\chi(G) \leq \aleph_0$, and thus by lemma 2.10, we get $\chi_{\boldsymbol{\Delta}_2^1}(G) = \chi_{\boldsymbol{\Sigma}_2^1}(G) = \chi(G)$. □

**Corollary 4.5.** *Let $S$ be a set of ordinals. In $\boldsymbol{L}[S]$, every locally countable Borel graph $G$ satisfies $\chi_{\boldsymbol{\Delta}_2^1}(G) = \chi(G)$.*

Since $S$ may be empty, this in particular holds in $\mathbf{L}$ itself.

We note that for $n = 2$, lemma 4.2 can be shown in a more concise way as follows. Given a vertex $v$, first modify a $\boldsymbol{\Delta}_2^1$ well-order as in the conclusion of lemma 4.3 to order proper colourings of the connected component of $v$ below improper colourings (using the enumeration of the connected component given by the well-order), and colourings using fewer colours below those using more. These properties of colourings are clearly $\boldsymbol{\Delta}^1_2$ because the graph is assumed to be Borel. Finally simply colour the vertex $v$ according to the colour it receives under the least colouring of its connected component using the modified well-order of $\mathcal{N}$.

While the results above concern graphs with few edges, note that the existence of a $\boldsymbol{\Delta}^1_2$ well-order of $\mathcal{N}$ implies the continuum hypothesis [Man] and hence that $\chi_{\boldsymbol{\Delta}^1_2}(G) = \chi(G)$ for any Polish graph $G$ of uncountable chromatic number.

In particular, for $\boldsymbol{\Delta}^1_2$ graphs in a model of ZFC with a $\boldsymbol{\Delta}^1_2$-definable well-order of the reals, we know the $\boldsymbol{\Delta}^1_2$ chromatic number unless our graph is *not* locally countable but has finite chromatic number.

**Problem 4.6.** *Let $G$ be a Borel graph that is not locally countable in a model of ZFC in which there is a $\boldsymbol{\Delta}^1_2$ well-order of $\mathcal{N}$ and such that $\chi(G) < \aleph_0$. Can it be true that $\chi(G) < \chi_{\boldsymbol{\Delta}^1_2}(G) < \chi_{\mathcal{B}}(G)$?*

The case $\chi(G) = \aleph_0$ is not relevant for this question because the continuum hypothesis would then imply $\chi_{\boldsymbol{\Delta}^1_2}(G) \in \{ \chi(G), \chi_{\mathcal{B}}(G) \}$.

For classes $\boldsymbol{\Delta}^1_3$ and above, the main obstacle to a result analogous to theorem 4.4 is the failure of lemma 4.3 since its proof relies on the correspondence between $\boldsymbol{\Sigma}^1_2$ and $\Sigma^{\mathrm{HC}}_1$. We recover a weaker result, namely that if there is a $\boldsymbol{\Delta}^1_n$ well-order of $\mathcal{N}$, then the $\boldsymbol{\Delta}^1_{n+1}$ chromatic number is equal to the classical chromatic number. We do this by once again applying lemma 4.2; the increase in complexity is due to the definition of the well-behaved well-order as follows.

**Lemma 4.7.** *Let $n \geq 2$, let $G$ be a $\boldsymbol{\Delta}^1_n$ graph, and assume there is a $\boldsymbol{\Delta}^1_n$-definable well-order of $\mathcal{N}$. Then there exists a $\boldsymbol{\Delta}^1_{n+1}$-definable well-order $\preceq$ of $\mathcal{N}$ with the property that $\preceq$ restricted to any one connected component of $G$ has order type at most $\omega$.*

*Proof.* Let $n \geq 2$, let $G = (V, E)$ be a $\boldsymbol{\Delta}^1_n$ graph, and let $\preceq$ be a $\boldsymbol{\Delta}^1_n$-definable well-order of $\mathcal{N}$.

We define a new well-order $\bar{\preceq}$ of $\mathcal{N}$ with the desired properties. Write $F := V^\omega$. By assumption, this is well-ordered by some $\boldsymbol{\Delta}^1_n$-definable $\preceq_F$.

Let $v \in V$. Define $e(v) = f$ if and only if $f \in F$ such that $f(\mathrm{ot}(\boldsymbol{\exists}_v)) = \boldsymbol{\exists}_v$, such that $f$ is injective, and such that $f$ is $\preceq_F$-minimal with these properties. Then $e : V \to F$ is $\boldsymbol{\Delta}^1_{n+1}$-definable, with the increase in complexity stemming from the minimality requirement, and for $v, w \in V$ with $w \in \boldsymbol{\exists}_v$, we have $e(w) = e(v)$. If $v, w$ are not connected, then we have $e(v)(0) \neq e(w)(0)$ and thus $e(v) \neq e(w)$.

We can now define for $v, w \in V$ that $v \mathrel{\bar{\prec}} w$ if $w \notin \boldsymbol{\exists}_v$ and $e(v) \prec_F e(w)$ or if $w \in \boldsymbol{\exists}_v$ and $\forall n, m \in \omega\,(e(v)(n) = v \wedge e(w)(m) = w \to n < m)$. By construction, this is a $\mathbf{\Delta}^1_{n+1}$-definable well-order with the desired property. $\square$

Note the consolation prize for the increase in complexity: the result works for any $\mathbf{\Delta}^1_n$ graph, whereas lemma 4.3 was only shown for Borel graphs.

**Theorem 4.8.** *Let $n \geq 2$. If there is a $\mathbf{\Delta}^1_n$ well-order of $\mathcal{N}$, then for every locally countable $\mathbf{\Delta}^1_n$ graph $G$, we have $\chi_{\mathbf{\Delta}^1_{n+1}}(G) = \chi(G)$.*

*Proof.* Let $n \geq 2$, and let $G$ be a locally countable $\mathbf{\Delta}^1_n$ graph. If there exists a $\mathbf{\Delta}^1_n$-definable well-order of $\mathcal{N}$, then by lemma 4.7 we can find a $\mathbf{\Delta}^1_{n+1}$-definable well-order $\preceq$ of $\mathcal{N}$ with the property that $\preceq$ restricted to any one connected component of $G$ has order type at most $\omega$. By lemma 4.2 and lemma 2.10, the claim follows. $\square$

## 4.1 Critical Pointclasses for Borel Graphs

Denote by $G_1$ the Borel graph which joins two elements $x, y \in 2^\omega$ precisely when there is $N < \omega$ such that $x \restriction (\omega \setminus N) = y \restriction (\omega \setminus N)$, where $f \restriction A = f \cap (A \times A)$ denotes the restriction of a function $f$ to a domain $A$. It is well known that this graph has uncountable Borel, and in fact also Baire-measurable and Lebesgue-measurable, chromatic number (folklore results following from the Kuratowski-Ulam theorem and Lebesgue density theorem, respectively). Since this graph is locally countable, by theorem 4.4, if there is a $\mathbf{\Delta}^1_2$ well-order of the reals then the $\mathbf{\Delta}^1_2$ chromatic number of $G_1$ is $\aleph_0$.[6]

**Theorem 4.9.** *There is a model of set theory in which $\chi_{\mathbf{\Delta}^1_1}(G_1) > \aleph_0$ but $\chi_{\mathbf{\Delta}^1_2}(G_1) = \aleph_0$.*

*Proof.* This follows directly from 4.4. $\square$

For higher pointclasses, the result is weaker.

**Theorem 4.10.** *Let $n \in \omega$, $n \geq 2$. If there is a model of set theory with $n-1$ Woodin cardinals and a measurable cardinal, then there is a model of set theory in which $\chi_{\mathbf{\Delta}^1_n}(G_1) > \aleph_0$ but $\chi_{\mathbf{\Delta}^1_{n+2}}(G_1) = \aleph_0$.*

*Proof.* For $n = 2$ we can actually weaken the assumptions: let $U$ be a normal ultrafilter witnessing that a cardinal $\kappa$ is measurable. Then in the model $\mathbf{L}[U]$ there is a $\boldsymbol{\Delta}^1_3$ well-order of the reals while every $\boldsymbol{\Sigma}^1_2$ set of reals (and hence also every $\boldsymbol{\Delta}^1_2$ set of reals) is Lebesgue-measurable (by a result of Solovay; see [KaA, theorem 14.2 and corollary 14.3]). Hence, by theorem 4.8, $\chi_{\boldsymbol{\Delta}^1_4}(G_1) = \aleph_0$, while $\chi_{\boldsymbol{\Delta}^1_2}(G_1) > \aleph_0$ in this model.

[6] It cannot be finite since $G_1$ trivially contains an $\aleph_0$-clique.

For $n > 2$ let $M$ be the canonical mouse with $n-1$ Woodin cardinals. Then $M$ satisfies that there is a $\boldsymbol{\Delta}^1_{n+1}$-definable well-order of the reals ([MaS]), but it in particular also has $n-2$ Woodin cardinals with a measurable above them, which implies that $\boldsymbol{\Pi}^1_{n-1}$ determinacy holds, whence every $\boldsymbol{\Pi}^1_n$ set of reals (and thus every $\boldsymbol{\Delta}^1_n$ set of reals) is Lebesgue-measurable (see, for example, [Mos, chapter 6]). Therefore in this model, $\chi_{\boldsymbol{\Delta}^1_{n+2}}(G_1) = \aleph_0$ while $\chi_{\boldsymbol{\Delta}^1_n}(G_1) > \aleph_0$. □

By the critical pointclass of a Borel graph $G$ of uncountable Borel chromatic number we mean the least $n$ such that the $\boldsymbol{\Delta}^1_n$ chromatic number of $G$ is countable, if such $n$ exists. In order to find critical pointclasses for Borel graphs in general, the gap in the theorem above must be closed.

**Problem 4.11.** *Let $n \geq 2$. If there exist $n-1$ Woodin cardinals, do we have $\chi_{\boldsymbol{\Delta}^1_{n+1}}(G_1) = \aleph_0$ or $\chi_{\boldsymbol{\Delta}^1_{n+1}}(G_1) > \aleph_0$?*

# 5 Open Questions

There is much to investigate in the context of projective colourings, and our investigation has left open far more problems than it has answered. We have mentioned problems 3.6 and 3.8 as well as the generalised conjecture 3.11. Furthermore, it would be interesting to determine whether there is a model where projective chromatic numbers are meaningfully distinct from all of their trivial bounds at the same time, a more general question than problem 4.6.

**Problem 5.1.** *Given $n \geq 2$, is there a model of set theory and a $\boldsymbol{\Delta}^1_n$ graph $G$ such that $\chi(G) < \chi_{\boldsymbol{\Delta}^1_n}(G) < \chi_{\mathcal{B}}(G)$? If there is, can we additionally require $\chi_{\boldsymbol{\Delta}^1_n}(G) \neq \chi_\mu(G)$? Furthermore, what is the consistency strength of the existence of such a graph?*

Further questions regard the strength of the theorems from section 4.

**Problem 5.2.** *Given $n \geq 2$, is there a model of set theory without a $\boldsymbol{\Delta}^1_n$ well-order of the reals such that for locally countable $\boldsymbol{\Delta}^1_n$ graphs, the $\boldsymbol{\Delta}^1_n$ chromatic number is equal to the classical chromatic number? If so, what is the consistency strength of this statement? Furthermore, can we extend the result to classes of graphs which are not necessarily locally countable?*

Note that theorem 4.8 holds for $\mathbf{\Delta}^1_n$ graphs, while the stronger theorem 4.4 concerns only Borel graphs.

**Problem 5.3.** *Let $G$ be a locally countable $\mathbf{\Delta}^1_2$ graph. If there is a $\mathbf{\Delta}^1_2$-definable well-order of $\mathcal{N}$, does this imply that $\chi_{\mathbf{\Delta}^1_2}(G) = \chi(G)$?*

Problem 4.11 concerns the critical pointclasses for projective chromatic numbers. In addition to this question, in the proof of theorem 4.10 the existence of finitely many Woodin cardinals is a very strong assumption. We conjecture that this assumption can be weakened considerably, in particular to the existence of an inaccessible cardinal.

**Problem 5.4.** *Given $n \geq 1$, what is the consistency strength of $\chi_{\mathbf{\Delta}^1_n}(G_1) > \aleph_0$ while $\chi_{\mathbf{\Delta}^1_{n+1}}(G_1) = \aleph_0$?*

## Acknowledgements

The first author would like to thank the members of the graduate school of system informatics at Kobe University, in particular Diego Mejia, for their valuable feedback. The second author is grateful to his advisors, Justin Moore and Paul Larson, for many insightful conversations about the problems addressed in this work. The work of the second author was supported by NSF grants DMS-2153975 and DMS-2451350.

# References

[Bar] T. Bartoszynski, H. Judah, *Set Theory: on the structure of the real line*, AK Peters/CRC Press, 1995.

[Ber] A. Bernshteyn, *Descriptive Combinatorics and Distributed Algorithms*, Notices of the American Mathematical Society 69.09 (2022).

[Con] C. T. Conley, A. S. Kechris, *Measurable chromatic and independence numbers for ergodic graphs and group actions*, Groups, Geometry, and Dynamics 7.1 (2013), 127–180.

[FFZ] V. Fischer, S. D. Friedman, L. Zdomskyy, *Cardinal characteristics, projective wellorders and large continuum*, Annals of Pure and Applied Logic 164.7-8 (2013), 763–770.

[Ges] S. Geschke, *Weak Borel Chromatic Numbers*, Mathematical Logic Quarterly 57(1) (2011), 5–13.

[Har] L. Harrington, *Long projective wellorderings*, Annals of Mathematical Logic 12.1 (1977).

[Jec] T. Jech, *Set Theory: The Third Millenium Edition*, Springer Berlin Heidelberg, 2003.

[KaA] A. Kanamori, *The Higher Infinite: Large Cardinals in Set Theory From Their Beginnings*, Springer Berlin Heidelberg, 2003.

[KaV] V. Kanovei, V. Lyubetsky, *Counterexamples to countable-section* $\Pi^1_2$ *uniformization and* $\Pi^1_3$ *separation*, Annals of Pure and Applied Logic 167.3 (2016), 262–283.

[Kec] A. Kechris, *Classical descriptive set theory*, Springer Science & Business Media, 2012.

[KST] A. Kechris, S. Solecki, S. Todorevic, *Borel Chromatic Numbers*, https://core.ac.uk/download/pdf/82238512.pdf (1999).

[Lyu] V. Lubetsky, *The existence of a nonmeasurable set of type* $A_2$ *implies the existence of an uncountable set of type* $CA$ *that does not contain a perfect subset*, Doklady Akademii Nauk SSSR 195 (1971), 548–550. Translated in Soviet Mathematics, Doklady 11 (1970), 1513–1515.

[Man] R. Mansfield, *The nonexistence of* $\mathbf{\Sigma}^1_2$ *wellorderings of the Cantor set*, Fundamenta Mathematicae 86 (1975), 279–282.

[MaS] D. A. Martin, J. R. Steel, *Iteration Trees*, J. Amer. Math. Soc. 7 (1994), no. 1, 1–73.

[Mos] Y. N. Moschovakis, *Descriptive set theory*, American Mathematical Society, 2025.

[Pal] K. Palamourdas, *1, 2, 3,..., $2n+1$, infinity!*, Diss. UCLA, 2012.

[Sol] R. M. Solovay, *The cardinality of* $\mathbf{\Sigma}^1_2$ *sets of reals*, In: J. J. Bulloff, T. C. Holyoke, and S. W. Hahn (eds.) *Foundations of Mathematics*. Symposium papers celebrating the sixtieth birthday of Kurt Gödel. Berlin, Springer-Verlag, 1969.